\documentclass[10pt]{article}
\usepackage{amsmath}
\usepackage{amssymb}
\usepackage{graphicx}
\newtheorem{theorem}{\scshape \mdseries  Theorem}[section]
\newtheorem{definition}[theorem]{\scshape \mdseries  Definition}
\newtheorem{lemma}[theorem]{\scshape \mdseries  Lemma}

\def \T{\mathcal{T}}
\def \A{\mathcal{A}}
\numberwithin{figure}{section}
\numberwithin{equation}{section}
\begin{document}

\title{\sf The linear unicyclic hypergraph with the second or third largest spectral radius\thanks{
This work was supported by Natural Science Foundation of China (11871073, 11871077), NSF of Department of Education of Anhui Province (KJ2017A362).}}
\author{Chao Ding$^{1,2}$,~
Yi-Zheng Fan$^{2,}$\thanks{Corresponding author.
 E-mail address: fanyz@ahu.edu.cn(Y.-Z. Fan), dcmath@sina.cn (C. Ding), 1500256209@qq.com (J.-C. Wan).
},~
Jiang-Chao Wan$^2$\\
{\small  \it $1$. School of Mathematics and Computational Science, Anqing Normal University, Anqing 246133, P. R. China} \\
 {\small  \it $2$. School of Mathematical Sciences, Anhui University, Hefei 230601, P. R. China}
}

\date{}
\maketitle

{\small
\noindent
\textbf{Abstract}:
The spectral radius of a uniform hypergraph is defined to be that of the adjacency tensor the hypergraph.
It is known that the unique unicyclic hypergraph with the largest spectral radius is a nonlinear hypergraph,
and the unique linear unicyclic hypergraph with the largest spectral radius is a power hypergraph.
In this paper we determine the the unique linear unicyclic hypergraph with the second or third largest spectral radius, where the former hypergraph is a power hypergraph and
the latter hypergraph is a non-power hypergraph.

%The spectral radius of hypergraph $G$ is the largest eigenvalue of its adjacency tensor, denoted by $\rho(G)$. If $G$ is not a power hypergraph obtained from a general graph, then $G$ is called a non g-power hypergraph. Let $K_{r,1}$ be a hyperstar with $r$ edges and a $r$-degree vertex $O$, $C_{i}$ be an unicycle with $i$ edges. Denote $C(i,r)$ be an unicycle hypergraph obtained from $C_{i}$ by attaching $K_{r,1}$ (at $O$ and a 1-degree vertex of $C_{i}$). Denote $C(i,T^\ast)$ be an unicycle hypergraph with $n$ edges obtained from $C_i$ by attaching some hypertrees at some 1-degree vertices of $C_i$. We prove that $\rho(C(i,T^\ast))\leq\rho(C(3,n-3)),$ the equality holds if only if $C(i,T^\ast)\cong{C(3,n-3)}.$ In this paper, we give the definition of subdividing an edge on hypergraphs. Let $G$ be obtained from $C(i,r)$ by subdividing a non-power edge $e$. We get that $\rho(G)>\rho(C(i,r))$ if $e$ lie on $K_{r,1}$, and $\rho(G)<\rho(C(i,r))$ if $e$ lie on $C_i$.


}

\section{Introduction}
A hypergraph $G=(V,E)$ consists of a nonempty vertex set $V=\{v_1,v_2,{\cdots},v_n\}$ denoted by $V(G)$ and a edge set $E=\{e_1,e_2,{\cdots},e_m\}$ denoted by $E(G)$,
 where $e_i \subseteq V$ for $i \in [m]:=\{1,2,\cdots,m\}$.
 If $|e_i|=k$ for each $i \in [m]$ and $k\geq2$, then $G$ is called a \emph{$k$-uniform hypergraph}.
 In particular, the $2$-uniform hypergraphs are exactly the classical simple graphs.
% A subgraph of hypergraph $G=(V,E)$ is $G'=(V',E')$ such that $V'{\subseteq}V$, $E'{\subseteq}E$.
The \emph{degree} of a vertex is the number of edges containing the vertex.
A vertex $v$ of $G$ is called a \emph{cored vertex} if it has degree one.
An edge $e$ of $G$ is called a \emph{pendent edge} if it contains $|e|-1$ cored vertices.
Sometimes a cored vertex in an pendent edge is also called a \emph{pendent vertex}.
A {\it walk} $W$ of length $l$ in $G$ is a sequence of alternate vertices and edges: $v_{0}e_{1}v_{1}e_{2}\ldots e_{l}v_{l}$,
    where $\{v_{i},v_{i+1}\}\subseteq e_{i}$ for $i=0,1,\ldots,l-1$.
If $v_0=v_l$, then $W$ is called a {\it circuit}.
A walk of $G$ is called a {\it path} if no vertices or edges are repeated.
A circuit $G$ is called a {\it cycle} if no vertices or edges are repeated except $v_0=v_l$.
The  hypergraph $G$ is said to be {\it connected} if every two vertices are connected by a walk.

 % A walk $W$ of length $k$ in $G$ is a sequence of alternate vertices and edges: $v_0e_1v_1e_2\cdots{e_kv_k}$, where $v_i, v_{i+1}\in{e_i}$ for $i=0,1,\cdots,k-1$. If no vertices or edges are repeated except $v_0=v_k$, the walk is called a cycle. $G$ is called unicyclic hypergraph if $G$ contains exactly one cycle.  For each $i{\in}V$, denote $E_i=\{e{\in}E:i{\in}e\}$, and degree $d_i=|E_i|$. If two vertices $i$ and $j$ are in the same edge $e$, then we denote $i{\sim}j$.

If $G$ is connected and acyclic, then $G$ is called a {\it hypertree} (also called supertree in \cite{LSQ} and other literatures).
It is known that a $k$-uniform hypertree on $n$ vertices has $\frac{n-1}{k-1}$ edges \cite[Proposition 4, p.392]{ber}.
If $G$ is connected and contains exactly one cycle, then $G$ is called a {\it unicyclic hypergraph}.
A $k$-uniform unicyclic hypergraph on $n$ vertices has $\frac{n}{k-1}$ edges \cite{FanTPL}.

Hu, Qi and Shao \cite{HQS} introduced a class of hypergraphs constructed from simple graphs.
Let $G=(V,E)$ be a simple graph. For any $k\geq3$, the {\it $k$-th power of $G$}, denoted by $G^{k}:=(V^{k},E^{k})$,
  is defined as the $k$-uniform hypergraph with the set of vertices $V^{k}:=V\cup{\{i_{e,1},\ldots,i_{e,k-2}|e\in E}\}$ and the set of edges
$E^{k}:={\{e\cup{{\{i_{e,1},\ldots,i_{e,k-2}}}\}|e\in E}\}$.
If a hypergraph can be obtained from the power of a simple graph, then we will such hypergraph a \emph{power hypergraph}.
  If for any two distinct edges $e_i,e_j$ of $G$, $|e_i\cap{e_j}|\leq1$, then $G$ is called a \emph{linear hypergraph}.
 It is known that all hypertrees and power hypergraphs are linear.
 For a unicyclic hypergraph $G$, if $G$ is linear, then the unique cycle of $G$ is a power of a cycle (as a simple graph) of length at least $3$;
 otherwise, $G$ contains a pair of edges sharing exactly two vertices which yields the unique cycle of $G$, and any other pair of edges shares at most one vertices.

The \emph{adjacency tensor} of a $k$-uniform hypergraph $G$ \cite{CD} on $n$ vertices is defined to be a $k$th order $n$ dimensional tensor $\mathcal{A}(G)=(a_{i_1i_2{\cdots}i_k})$, where
$$a_{i_1i_2{\cdots}i_k}= \left  \{
\begin{array}{ll}
\frac{1}{(k-1)!}, & {\rm if}\ \{i_1,i_2,{\cdots},i_k\}{\in}E(G)  ,\\
0, & {\rm orthwise}.
\end{array}
\right.
$$
Qi \cite{qi05} introduces the eigenvalues of a supersymmetric tensor, from which one can get the definition of the eigenvalues of the adjacency tensor of a uniform hypergraph.
The \emph{spectral radius} of a uniform hypergraph is the maximum modulus of the eigenvalues of its adjacency tensor; see  more in Section 2.
%The adjacency tensor $\mathcal{A}(G)$ is nonnegative, and it is irreducible if and only if $G$ is connected. In this paper, we only consider connect linear uniform hypergraph and denote $\rho(G)=\rho(\mathcal{A}(G))$.

The spectral hypergraph theory has emerged as a hot topic in algebraic graph theory \cite{BL,CD,FHB, FBH, LM,PZ,XW19, ZKSB,zhou}.
Among all uniform hypertrees with given number of vertices or edges, researchers worked on the ordering the hypertrees by their spectral radii.
In 2015, Li, Shao and Qi \cite{LSQ} determined the hypertrees with the largest and the second largest spectral radii respectively.
In 2016, Yuan, Shao and Shan \cite{YSS} determined the first eight hypertrees with largest spectral radii,
and in 2017 Yuan, Si and Zhang \cite{YSZ} determined the ninth and tenth hypertrees with largest spectral radii.
In 2016, Fan, Tan, Peng and Liu \cite{FanTPL} investigated the hypergraphs that attain largest spectral radii among all hypergraphs with given number of edges.
They determined the unique unicyclic hypergraphs with the largest spectral radius, which is not a linear hypergraph;
and they also determined the unique linear unicyclic  hypergraph with the largest spectral radius, which is a power hypergraph.
They proposed several candidates for the linear  bicyclic hypergraph with the largest spectral radius.
Later in 2018, Kang et al. \cite{KLQY} confirmed a conjecture in \cite{FanTPL} which lead to the unique linear  bicyclic hypergraph with the largest spectral radius.
Recently, Ouyang, Qi and Yuan \cite{OQY} considered the nonlinear hypergraphs, and determined the first five unicyclic hypergraphs and first three bicyclic hypergraphs with largest spectral radii.
Other works on the ordering of hypertrees or unicylic hypergraphs can be referred  to \cite{CCZ, GZ, SKLS18, SKLS,XW, XWD, ZC}.

In this paper we continue the work on the ordering of linear unicyclic hypergraphs by their spectral radii, and
determine the the unique linear unicyclic hypergraph with the second or third largest spectral radius, where the former hypergraph is a power hypergraph and
the latter hypergraph is a non-power hypergraph.

%
%Recently, the research on spectrum radius of hypergraphs has attracted great interests. By using the operation of moving edge, Li, Shao and Qi {\cite{1}} proved that the hypertree with the maximal spectral radius is the hyperstar. Lu and Man {\cite{2}} designed the technique of weighted incidence matrix and showed that the smallest limit point of the spectral radius of connected $m$-uniform hypergraphs is $\sqrt[m]{4} $. Kang et al. {\cite{3}} determined the largest spectral radius of linear bicyclic uniform hypergraths. Zhang et al. {\cite{4}} extended the technique of weighted incidence matrix and determined the hypergraphs with the maximum spectral radius over all nearly uniform supertrees. There are more results in {\cite{5}-\cite{20}}. By using the operation of moving edge, Fan et al.{\cite{2}} showed that the unicycle linear hypergraph with the maximum spectral radius is $C'(3,n-3)$ (an unicycle hypergraph obtained from $C_{3}$ by attaching $K_{n-3,1}$ (at the center vertex of $K_{n-3,1}$ and a 2-degree vertex of $C_{3}$), where $C'(3,n-3)$ is a g-power hypergraphs.  In this paper we describe the upper bound about the spectral radius of a kind of unicyclic non g-power hypergraphs by using the operation of moving edge and weighted incidence matrix. We give the definition of subdividing an edge on hypergraphs. We characterize the influence of subdividing edge on the spectral radius of a kind of unicyclic hypergraphs.

\section{Preliminaries}

For integers $k\geq 3$ and $n\geq 2$,
  a real {\it tensor} (also called {\it hypermatrix}) $\mathcal{T}=(t_{i_{1}\ldots i_{k}})$ of order $k$ and dimension $n$ refers to a
  multidimensional array with entries $t_{i_{1}i_2\ldots i_{k}}$ such that $t_{i_{1}i_2\ldots i_{k}}\in \mathbb{R}$ for all $i_{j}\in [n]$ and $j\in [k]$.
 The tensor $\mathcal{T}$ is called \textit{symmetric} if its entries are invariant under any permutation of their indices.
 Given a vector $x\in \mathbb{R}^{n}$, $\mathcal{T}x^{k}$ is a real number, and $\mathcal{T}x^{k-1}$ is an $n$-dimensional vector, which are defined as follows:
   $$\mathcal{T}x^{k}=\sum_{i_1,i_{2},\ldots,i_{k}\in [n]}t_{i_1i_{2}\ldots i_{k}}x_{i_1}x_{i_{2}}\cdots x_{i_k},$$
   $$   (\mathcal{T}x^{k-1})_i=\sum_{i_{2},\ldots,i_{k}\in [n]}t_{ii_{2}i_3\ldots i_{k}}x_{i_{2}}x_{i_3}\cdots x_{i_k}, \mbox{~for~} i \in [n].$$
 Let $\mathcal{I}$ be the {\it identity tensor} of order $k$ and dimension $n$, that is, $i_{i_{1}i_2 \ldots i_{k}}=1$ if and only if
   $i_{1}=i_2=\cdots=i_{k} \in [n]$ and zero otherwise.

\begin{definition}{\em\cite{qi05,CPZ2}} Let $\mathcal{T}$ be a $k$-th order $n$-dimensional real tensor.
For some $\lambda \in \mathbb{C}$, if the polynomial system $(\lambda \mathcal{I}-\mathcal{T})x^{k-1}=0$, or equivalently $\mathcal{T}x^{k-1}=\lambda x^{[k-1]}$, has a solution $x\in \mathbb{C}^{n}\backslash \{0\}$,
then $\lambda $ is called an eigenvalue of $\mathcal{T}$ and $x$ is an eigenvector of $\mathcal{T}$ associated with $\lambda$,
where $x^{[k-1]}:=(x_1^{k-1}, x_2^{k-1},\ldots,x_n^{k-1}) \in \mathbb{C}^n$.
\end{definition}

If $x$ is a real eigenvector of $\mathcal{T}$, surely the corresponding eigenvalue $\lambda$ is real.
In this case, $x$ is called an {\it $H$-eigenvector} and $\lambda$ is called an {\it $H$-eigenvalue}.
%Furthermore, if $x\in \mathbb{R}_{+}^{n}$ (the set of nonnegative vectors of dimension $n$), then $\lambda $ is called an {\it $H^{+}$-eigenvalue} of $\mathcal{T}$;
%if $x\in \mathbb{R}_{++}^{n}$ (the set of positive vectors of dimension $n$), then $\lambda$ is said to be an {\it $H^{++}$-eigenvalue} of $\mathcal{T}$.
The {\it spectral radius of $\T$} is defined as
$$\rho(\T)=\max\{|\lambda|: \lambda \mbox{ is an eigenvalue of } \T \}.$$

By the Perron-Frobenius theorem for nonnegative tensors\cite{CPZ, FGH, YY}, the spectral radius of $\A(G)$, also referred to the \emph{spectral radius of $G$}, denoted by $\rho(G)$, is exactly the largest $H$-eigenvalue of $\A(G)$.
If $G$ is connected, there exists a unique positive eigenvector up to scales corresponding to $\rho(G)$, called the {\it Perron vector} of $G$.
%In addition, $\rho(G)$ is the optimal value of the following maximization (see \cite{PZ}):
%$$ \rho(G)=\max_{x \in \mathbb{R}^n, \|x\|_k=1}\A(G)x^k=\max_{x \in \mathbb{R}^n, \|x\|_k=1}\sum_{e=\{u_1,u_2,\cdots, u_k\} \in E(G)}k x_{u_1}x_{u_2}\cdots x_{u_k}.\eqno(2.1)$$
%The eigenvector equation  $\mathcal{A}(G)x^{k-1}=\lambda x^{[k-1]}$ could be interpreted as
%$$ \lambda x_u^{k-1}= \sum_{\{u,u_2,u_3,\ldots, u_k\} \in E(G)} x_{u_2}x_{u_3} \cdots x_{u_k}, \mbox{~for each~} u \in V(G).\eqno(2.2)$$

Li, Shao and Qi \cite{LSQ} introduce the operation of {\it moving edges} on hypergraphs.
Let $r \ge  1$ and let $G$ be a hypergraph with $u \in V(G)$ and $e_1, \ldots, e_r \in E(G)$ such that $u \notin e_i$ for $i = 1, \ldots, r$.
Suppose that $v_i \in e_i$ and write $e'_i = (e_i \backslash \{v_i\}) \cup \{u\}$ ($i = 1, \ldots, r)$.
Let $G'$ be the hypergraph with $V(G')=V(G)$ and $E(G') = (E \backslash \{e_1, \ldots, e_r\}) \cup \{e'_1, \ldots, e'_r\}$.
We say that $G'$ is obtained from $G$ by moving edges $(e_1, \ldots, e_r)$ from $(v_1, \ldots, v_r)$ to $u$.

\begin{lemma} \label{mov} {\em \cite{LSQ}}
Let $r \ge  1$ and let $G$ be a connected hypergraph.
 Let $G'$ be obtained from $G$ by moving edges $(e_1, \ldots, e_r)$ from $(v_1, \ldots, v_r)$ to $u$.
 Assume that $G'$ contains no multiple edges.
If $x$ is a Perron vector of $G$ and $x_u \ge \max_{1\le i \le r}x_{v_i}$, then $\rho(G') > \rho(G)$.
\end{lemma}

\begin{figure}
\centering
  \setlength{\unitlength}{1bp}
  \begin{picture}(411.70, 125.21)(0,0)
  \put(0,0){\includegraphics{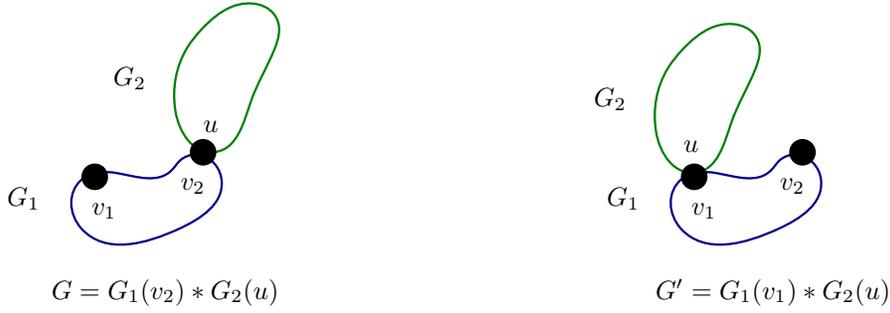}}
  \put(37.50,39.63){\fontsize{10}{12}\selectfont $v_1$}
  \put(70.97,49.11){\fontsize{10}{12}\selectfont $v_2$}
  \put(79.45,71.04){\fontsize{10}{12}\selectfont $u$}
  \put(45.55,88.49){\fontsize{10}{12}\selectfont $G_2$}
  \put(5.67,43.12){\fontsize{10}{12}\selectfont $G_1$}
  \put(263.47,39.63){\fontsize{10}{12}\selectfont $v_1$}
  \put(296.87,49.11){\fontsize{10}{12}\selectfont $v_2$}
  \put(260.57,63.29){\fontsize{10}{12}\selectfont $u$}
  \put(226.68,80.73){\fontsize{10}{12}\selectfont $G_2$}
  \put(231.56,43.12){\fontsize{10}{12}\selectfont $G_1$}
  \put(22.29,8.12){\fontsize{10}{12}\selectfont $G=G_{1}(v_{2})\ast G_{2}(u)$}
  \put(250.00,8.12){\fontsize{10}{12}\selectfont $G'=G_{1}(v_{1})\ast G_{2}(u)$}
  \end{picture}
\caption{An illustration of relocating subhypergraph}\label{move}
\end{figure}

Fan et al. \cite{FanTPL} introduced a special case of moving edges.
Let $G_1,G_2$ be two vertex-disjoint hypergraphs, where $v_{1},v_{2}$ are two distinct vertices of $G_{1}$ and $u$ is a vertex of $G_{2}$ (called the root of $G_2$).
Let $G=G_{1}(v_{2})\ast G_{2}(u)$ (respectively, $G'=G_{1}(v_{1})\ast G_{2}(u)$) be the hypergraph obtained by identifying $v_2$ with $u$ (respectively, identifying $v_1$ with $u$);
see the hypergraphs in Fig. \ref{move}.
It is said that $G'$ is obtained from $G$ by {\it relocating} $G_{2}$ rooted at $u$ from $v_{2}$ to $v_{1}$.

\begin{lemma}\label{relocate} {\em \cite{FanTPL}}
Let $G=G_{1}(v_{2})\ast G_{2}(u)$ and $G'=G_{1}(v_{1})\ast G_{2}(u)$ be two connected hypergraphs.
If there exists a Perron vector $x$ of $G$ such that $x_{v_1} \ge x_{v_2}$, then $\rho(G') > \rho(G)$.
\end{lemma}

Yuan, Shao and Shan \cite{YSS} defined an new type of edge-moving operation.

\begin{definition}\label{newmov}{\em \cite{YSS}}
Let $e,f$ be two edges of a $k$-uniform connected hypergraph $G$ such that $e \cap f =V_1$, where $|V_1|=k-r$, $2 \le r \le k-1$.
Write $e \backslash V_1=\{u_1, \ldots, u_r\}$ and  $f \backslash V_1=\{v_1, \ldots, v_r\}$,
where $u_1$ and $v_1$ are non-pendent vertices, but $u_2, \ldots, u_r$ and $v_2, \ldots, v_r$ are all pendent vertices.
Define $G_{e,f}$ be the hypergraph obtained from $G$ by moving all edges incident to $v_1$ except $f$ from $v_1$ to $u_2$.
\end{definition}

\begin{lemma}{\em \cite{YSS}}\label{Ymov}
Let $G$ be a $k$-uniform connected hypergraph, and let $e, f$ be two edges of $G$ satisfying the condition in Definition \ref{newmov}.
Then $\rho(G_{e,f})>\rho(G)$.
\end{lemma}

Lu and Man \cite{LM} introduced a novel method for computing or comparing the spectral radii of hypergraphs.

\begin{definition} {\em \cite{LM}}
A weighted incidence matrix $B$ of a hypergraph $G=(V,E)$ is a $|V|\times|E|$ matrix such that for any vertex $v$ and any edge $e$, the entry $B(v,e)>0$ if $v\in{e}$ and $B(v,e)=0$ if $v\notin{e}$.
\end{definition}

\begin{definition} {\em \cite{LM}} Let $G$ be hypergraph with a weighted incidence matrix $B$.

$(1)$ $G$ is called $\alpha$-normal if $B$ satisfies

 \mbox{\em(i)} $\sum_{e: v \in e}B(v,e)=1$, for any $v \in V(G)$,

\mbox{\em(ii)} $\prod_{v \in e}B(v, e)=\alpha$, for any $e \in E(G)$.

%$(2)$ $G$ is called $\alpha$-subnormal if $B$ satisfies
%
%\mbox{\em(i)}  $\sum_{e: v \in e}B(v,e) \le 1$, for any $v \in V(G)$.
%
%\mbox{\em(ii)} $\prod_{v \in e}B(v, e) \ge \alpha$, for any $e \in E(G)$.
%
%Moreover $G$ is called strictly $\alpha$-subnormal if $G$ is $\alpha$-subnormal but not $\alpha$-normal.

$(2)$ $G$ is called $\alpha$-supernormal if $B$ satisfies

\mbox{\em(i)}  $\sum_{e: v \in e}B(v,e) \ge 1$, for any $v \in V(G)$.

\mbox{\em(ii)}  $\prod_{v \in e}B(v, e) \le \alpha$, for any $e \in E(G)$.

Moreover $G$ is called strictly $\alpha$-supernormal if $G$ is $\alpha$-supernormal but not $\alpha$-normal.

$(3)$ The incidence matrix $B$ is called consistent if for any cycle $v_0e_1v_1e_2\cdots e_l{v}_l\;(v_l=v_0)$
$$\prod\limits_{i=1}^l\frac{B(v_i,e_i)}{B(v_{i-1},e_i)}=1.$$
In this case, we call $G$ consistently $\alpha$-normal (resp. consistently  $\alpha$-supernormal) if $G$ is also  $\alpha$-normal (resp. $\alpha$-supernormal).
\end{definition}

\begin{lemma} {\em \cite{LM}} \label{comp}
 Let $G$ be a connected $k$-uniform hypergraph. Then the following results hold.

\mbox{\em(1)}
 $G$ is consistently $\alpha$-normal if and only if $\rho(G)=\alpha^{-1/k}$.

%\mbox{\em(2)}
% If $G$ is strictly $\alpha$-subnormal, then $\rho(G)<\alpha^{-1/k}$.

\mbox{\em(2)}
 If $G$ is strictly and consistently  $\alpha$-supernormal, then $\rho(G)>\alpha^{-1/k}$.

\end{lemma}

\begin{lemma}{\em \cite{zhou}} \label{pow}
Let $G^k$ be the $k$-th power of a simple graph $G$.
Then $\rho(G^k)=\rho(G)^{\frac{2}{k}}$.
\end{lemma}

%\begin{lemma} {\em \cite{10}} \label{pen10}
% If $G$ is a hypergraph with the maximum spectral radius among the connected general hypergraphs with fixed number of edges, then $G$ contains a vertex adjacent to all the other vertices.
%\end{lemma}

\section{Main results}

We first introduce some special graphs and hypergraphs.
Let $K_{1,s}$ be a star on $1+s$ vertices, and let $C_n$ be cycle of length $n$, both as simple graphs.
Let $S_{m,g}$ be a unicyclic graph obtained from a cycle $C_g$ by attaching a star $K_{1,m-g}$ at some vertex.
Let $T_{m,1}$ (respectively, $T_{m,2}$) be obtained from $S_{m-1,3}$
  (respectively, $S_{m-2,3}$) by attaching one pendent edge (respectively, two pendent edges) at some vertex of degree $2$.
Let $U_{m,1}$ be be obtained from $S_{m-1,3}$ by attaching one pendent edge at some pendent vertex.

The $k$-th power $K_{1,s}^k$ of $K_{1,s}$, is called a \emph{hyperstar} with $s$ edges, where the vertex of maximum degree is called the \emph{center} of the hyperstar.
Let $O_{m}$ be the $k$-uniform hypergraph obtained from the power $C^k_3$ by attaching a hyperstar $K_{1,m-3}^k$ with its center at some cored vertex.
Let $Q_{m}$ (respectively, $P_{m}$) be the $k$-uniform hypergraph obtained from the power $S^k_{m-1,3}$ by attaching a pendent edge
    to a cored vertex on the cycle adjacent to (respectively, not adjacent to) the vertex with maximum degree of $S^k_{m-1,3}$.

\begin{figure}
\centering
  \setlength{\unitlength}{1bp}%
   \begin{picture}(364.97, 220.07)(0,0)
  \put(0,0){\includegraphics{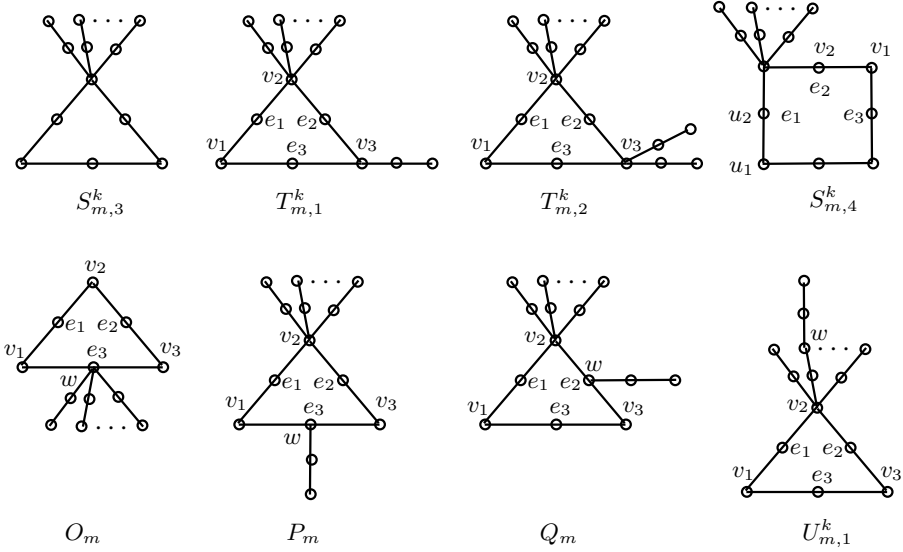}}
  \put(38.66,200.86){\fontsize{11.38}{13.66}\selectfont $\cdots$}
  \put(113.99,200.86){\fontsize{11.38}{13.66}\selectfont $\cdots$}
  \put(120.75,102.44){\fontsize{11.38}{13.66}\selectfont $\cdots$}
  \put(213.43,102.44){\fontsize{11.38}{13.66}\selectfont $\cdots$}
  \put(39.88,47.49){\fontsize{11.38}{13.66}\selectfont $\cdots$}
  \put(292.16,205.51){\fontsize{11.38}{13.66}\selectfont $\cdots$}
  \put(32.92,132.44){\fontsize{9.10}{10.93}\selectfont $S^k_{m,3}$}
  \put(108.25,132.44){\fontsize{9.10}{10.93}\selectfont $T^k_{m,1}$}
  \put(309.30,133.49){\fontsize{9.10}{10.93}\selectfont $S^k_{m,4}$}
  \put(112.27,7.77){\fontsize{9.10}{10.93}\selectfont $P_{m}$}
  \put(207.69,7.81){\fontsize{9.10}{10.93}\selectfont $Q_{m}$}
  \put(28.67,7.63){\fontsize{9.10}{10.93}\selectfont $O_{m}$}
  \put(298.42,165.81){\fontsize{9.10}{10.93}\selectfont $e_1$}
  \put(309.22,176.42){\fontsize{9.10}{10.93}\selectfont $e_2$}
  \put(279.08,146.52){\fontsize{9.10}{10.93}\selectfont $u_1$}
  \put(278.67,165.83){\fontsize{9.10}{10.93}\selectfont $u_2$}
  \put(332.13,191.38){\fontsize{9.10}{10.93}\selectfont $v_1$}
  \put(321.96,166.05){\fontsize{9.10}{10.93}\selectfont $e_3$}
  \put(310.91,191.38){\fontsize{9.10}{10.93}\selectfont $v_2$}
  \put(28.70,87.19){\fontsize{9.10}{10.93}\selectfont $e_1$}
  \put(40.76,87.19){\fontsize{9.10}{10.93}\selectfont $e_2$}
  \put(36.90,77.30){\fontsize{9.10}{10.93}\selectfont $e_3$}
  \put(5.67,77.47){\fontsize{9.10}{10.93}\selectfont $v_1$}
  \put(36.18,109.38){\fontsize{9.10}{10.93}\selectfont $v_2$}
  \put(64.47,78.25){\fontsize{9.10}{10.93}\selectfont $v_3$}
  \put(103.91,163.88){\fontsize{9.10}{10.93}\selectfont $e_1$}
  \put(115.97,163.88){\fontsize{9.10}{10.93}\selectfont $e_2$}
  \put(112.11,153.99){\fontsize{9.10}{10.93}\selectfont $e_3$}
  \put(102.95,180.76){\fontsize{9.10}{10.93}\selectfont $v_2$}
  \put(81.93,155.41){\fontsize{9.10}{10.93}\selectfont $v_1$}
  \put(138.41,155.44){\fontsize{9.10}{10.93}\selectfont $v_3$}
  \put(203.11,65.25){\fontsize{9.10}{10.93}\selectfont $e_1$}
  \put(215.16,65.25){\fontsize{9.10}{10.93}\selectfont $e_2$}
  \put(211.30,55.84){\fontsize{9.10}{10.93}\selectfont $e_3$}
  \put(201.90,82.37){\fontsize{9.10}{10.93}\selectfont $v_2$}
  \put(180.31,56.00){\fontsize{9.10}{10.93}\selectfont $v_1$}
  \put(238.83,56.00){\fontsize{9.10}{10.93}\selectfont $v_3$}
  \put(109.57,82.37){\fontsize{9.10}{10.93}\selectfont $v_2$}
  \put(89.06,57.84){\fontsize{9.10}{10.93}\selectfont $v_1$}
  \put(145.82,57.86){\fontsize{9.10}{10.93}\selectfont $v_3$}
  \put(110.29,65.25){\fontsize{9.10}{10.93}\selectfont $e_1$}
  \put(122.35,65.25){\fontsize{9.10}{10.93}\selectfont $e_2$}
  \put(118.49,55.84){\fontsize{9.10}{10.93}\selectfont $e_3$}
  \put(26.77,65.73){\fontsize{9.10}{10.93}\selectfont $w$}
  \put(111.19,43.89){\fontsize{9.10}{10.93}\selectfont $w$}
  \put(224.72,71.55){\fontsize{9.10}{10.93}\selectfont $w$}
  \put(312.04,77.02){\fontsize{11.38}{13.66}\selectfont $\cdots$}
  \put(306.56,7.82){\fontsize{9.10}{10.93}\selectfont $U^k_{m,1}$}
  \put(301.96,40.05){\fontsize{9.10}{10.93}\selectfont $e_1$}
  \put(314.02,40.05){\fontsize{9.10}{10.93}\selectfont $e_2$}
  \put(310.16,30.16){\fontsize{9.10}{10.93}\selectfont $e_3$}
  \put(300.76,56.93){\fontsize{9.10}{10.93}\selectfont $v_2$}
  \put(335.98,31.60){\fontsize{9.10}{10.93}\selectfont $v_3$}
  \put(280.22,31.57){\fontsize{9.10}{10.93}\selectfont $v_1$}
  \put(309.61,83.42){\fontsize{9.10}{10.93}\selectfont $w$}
  \put(213.69,200.86){\fontsize{11.38}{13.66}\selectfont $\cdots$}
  \put(207.95,132.44){\fontsize{9.10}{10.93}\selectfont $T^k_{m,2}$}
  \put(203.61,163.88){\fontsize{9.10}{10.93}\selectfont $e_1$}
  \put(215.67,163.88){\fontsize{9.10}{10.93}\selectfont $e_2$}
  \put(211.81,153.99){\fontsize{9.10}{10.93}\selectfont $e_3$}
  \put(202.65,180.76){\fontsize{9.10}{10.93}\selectfont $v_2$}
  \put(181.63,155.41){\fontsize{9.10}{10.93}\selectfont $v_1$}
  \put(238.11,155.44){\fontsize{9.10}{10.93}\selectfont $v_3$}
  \end{picture}%
  \caption{Eight linear unicyclic $k$-uniform hypergraphs with $m$ edges}\label{8graphs}
\end{figure}

\begin{lemma}\label{QT}
For $m \ge 4$, $\rho(Q_{m})< \rho(T^k_{m,1})$.
\end{lemma}

\emph{Proof.} Label the partial vertices and edges of $Q_m$ as in Fig. \ref{8graphs},
   where $v_1e_1v_2e_2v_3e_3v_1$ is the $3$-cycle, and $w$ is the vertex of the edge $e_2$ to which a pendent edge is attached.
Let $x$ be a Perron vector of $Q_{m}$.
If $x_{v_3} \ge x_w$, moving the pendent edge attached at $w$ from $w$ to $v_3$, we arrive at the hypergraph $T^k_{m,1}$.
By Lemma \ref{mov}, $\rho(Q_{m})< \rho(T^k_{m,1})$.
Otherwise, $x_w > x_{v_3}$, moving the edge $e_3$ from $v_3$ to $w$, we also arrive at the hypergraph $T^k_{m,1}$.
So, By Lemma \ref{mov}, $\rho(Q_{m})< \rho(T^k_{m,1})$.
The result follows. \hfill $\square$

\begin{lemma}\label{PQ}
For $m \ge 5$, $\rho(P_{m})<\rho(Q_{m})$.
\end{lemma}

\emph{Proof.} Label the partial vertices and edges of $P_m$ as in Fig. \ref{8graphs},
   where $v_1e_1v_2e_2v_3e_3v_1$ is the $3$-cycle, and $w$ is the vertex of the edge $e_3$ to which a pendent edge is attached.
We first construct a consistently $\alpha$-normal weighted incidence matrix $B$ of $P_{m}$.
Let $r:=m-4 \ge 1$, the number of pendent edges attached at $v_2$.
For each cored vertex $v$ incident to the unique edge $e$, $B(v,e)=1$.
For each pendent edge $e$ attached at $v \in \{v_2,w\}$), define $B(v,e)=\alpha$.
Define
\begin{align*}
B(w,e_3)&=1-\alpha, B(v_1,e_3)=B(v_3,e_3)=\sqrt{\frac{\alpha}{1-\alpha}}=:\beta,\\
B(v_1,e_1)&=B(v_3,e_2)=1-\beta, B(v_2,e_1)=B(v_2,e_2)=\frac{\alpha}{1-\beta}.
\end{align*}
Then $B$ is consistent.
Let
\begin{equation}\label{P}
f_P(\alpha):=\frac{2 \alpha}{1-\sqrt{\frac{\alpha}{1-\alpha}}}+r \alpha.
\end{equation}
Then $P_{m}$ is consistently $\alpha_0$-normal if $f_P(\alpha)=1$ has a solution $\alpha_0 \in (0,\frac{1}{2})$.
Observe that $f_P(\alpha) \to 0+$ if $\alpha \to 0+$, $f_P(\alpha) \to +\infty$ if $\alpha \to \frac{1}{2}-$, and
$f_P(\alpha)$ is strictly increasing in $(0,\frac{1}{2})$.
So $f_P(\alpha)=1$ has a unique solution $\alpha_0 \in (0,\frac{1}{2})$, and $\rho(P_{m})=\alpha_0^{-\frac{1}{k}}$ by Lemma \ref{comp}.
As $f_P(\frac{1}{5}) \ge 1$ and $f_P(\frac{1}{r+2})>1$,
\begin{equation} \label{upP}
\alpha_0 \le \frac{1}{5}, \alpha_0 < \frac{1}{r+2}.
\end{equation}

We next define a weighted incident matrix $\bar{B}$ of $Q_{m}$.
%Let $\alpha=\alpha_0$ in the following.
For each cored vertex $v$ incident to the unique edge $e$, $\bar{B}(v,e)=1$.
For each pendent edge $e$ attached at $v \in \{v_2,w\}$, define $\bar{B}(v,e)=\alpha$.
Define
\begin{align*}
\bar{B}(v_3,e_2)&=x, \bar{B}(v_2,e_2)=\frac{\alpha}{(1-\alpha)x}, \bar{B}(v_3,e_3)=1-x, \\
\bar{B}(v_1,e_3)&=\frac{\alpha}{1-x}=:\beta, \bar{B}(v_1,e_1)=1-\beta, \bar{B}(v_2,e_1)=\frac{\alpha}{1-\beta}.
\end{align*}
To make $\bar{B}$ be strictly consistently $\alpha$-supernormal, we need
\begin{equation}\label{inci-consi}
(1-\beta) \cdot \frac{\alpha}{(1-\alpha)x} \cdot (1-x)=\beta  \cdot \frac{\alpha}{1-\beta} \cdot  x.
\end{equation}
\begin{equation} \label{inci-v}
 h(x):=\frac{\alpha}{1-\beta}+\frac{\alpha}{(1-\alpha)x}+r \alpha >1,
\end{equation}
By Eq. (\ref{inci-consi}), we have
%\begin{equation}
%1-\beta=\frac{x}{1-x}\sqrt{\alpha(1-\alpha)}.
%\end{equation}
%So,
\begin{equation}\label{x}
x=\frac{1-\alpha}{1+\sqrt{\alpha(1-\alpha)}}=:\gamma.
\end{equation}
Now substituting (\ref{x}) to $h(x)$ by taking $\alpha=\alpha_0$, and combining Eq. (\ref{P}) and the fact $f_P(\alpha_0)=1$, we have
\begin{align*}
h(\gamma) &=  \frac{2\alpha_0}{1-\alpha_0}\sqrt{\frac{\alpha_0}{1-\alpha_0}}+\frac{\alpha_0 (2-\alpha_0)}{(1-\alpha_0)^2}+r \alpha_0 \\
&= 1+ \frac{2\alpha_0}{1-\alpha_0}\sqrt{\frac{\alpha_0}{1-\alpha_0}}+\frac{\alpha_0 (2-\alpha_0)}{(1-\alpha_0)^2}-\frac{2 \alpha_0}{1-\sqrt{\frac{\alpha_0}{1-\alpha_0}}}\\
& = 1+\frac{\alpha_0^2}{(1-\alpha_0)^2(1-2\alpha_0)}\left(1-4 \alpha_0 +2 \alpha_0^2 - 2 \alpha_0 \sqrt{\alpha_0(1-\alpha_0)}\right)\\
& =1+\frac{\alpha_0^2}{(1-\alpha_0)^2(1-2\alpha_0)}\left(-1+2\sqrt{1-\alpha_0}\left((1-\alpha_0)^{3 \over 2}-\alpha_0^{3 \over 2}\right)\right).
\end{align*}
Let $\phi(\alpha):=-1+2\sqrt{1-\alpha}\left((1-\alpha)^{3 \over 2}-\alpha^{3 \over 2}\right)$.
As $\alpha_0 \le \frac{1}{5}$ by (\ref{upP}), $\phi(\alpha_0) \ge \phi(\frac{1}{5})=\frac{3}{25}>0$, and hence $h(\gamma)>1$.
%When $r=1$, by Eq. (\ref{O}) we get $\alpha_0=\frac{1}{5}$, and $\phi(\frac{1}{5})=\frac{3}{25}>0$.
%When $r=2$, we get $h(\frac{1}{6})<1$ and hence $\alpha_0<\frac{1}{6}$, and then $\phi(\alpha_0)>\phi(\frac{1}{6})>0$.
%When $r \ge 3$, then $\alpha_0 <\frac{1}{r+2}\le \frac{1}{5}$, and surely $\phi(\alpha_0)>\phi(\frac{1}{5})>0$.
So, $Q_{m}$ is strictly consistently $\alpha_0$-supernormal, and by Lemma \ref{comp}
$$\rho(Q_{m})>\alpha_0^{-\frac{1}{k}}=\rho(P_{m}).$$
The result follows.\hfill $\square$

\begin{lemma}\label{OP}
For $m \ge 5$, $\rho(O_{m})< \rho(P_{m})$.
\end{lemma}

\emph{Proof.}
Label the partial vertices and edges of $O_m$ as in Fig. \ref{8graphs},
   where $v_1e_1v_2e_2v_3e_3v_1$ is the $3$-cycle, and $w$ is the vertex of the edge $e_3$ to which a hyperstar is attached.
We define a weighted incidence matrix $B$ of $O_m$ as follows.
For each pendent vertex $v$ incident to the unique edge $e$, $B(v,e)=1$.
%Let $w$ be the vertex to which all pendent edges are attached.
For each pendent edge $e$ attached at $w$, define $B(w,e)=\alpha$.
Let $r:=m-4 \ge 1$.
Define
\begin{align*}
B(w,e_3)&=1-(r+1)\alpha, B(v_1,e_3)=B(v_3,e_3)=\beta, \\
B(v_1,e_1)& =B(v_3,e_2)=1-\beta, B(v_2,e_1)=B(v_2,e_2)=\frac{\alpha}{1-\beta}.
\end{align*}
It is easily seen $B$ is consistent.
To make $O_{m}$ be $\alpha$-normal, we require
$$\beta=\sqrt{\frac{\alpha}{1-(r+1)\alpha}}, \; \frac{2\alpha}{1-\beta}=1.$$
%by considering (1)(ii) on the edge $e_3$ and (1)(i) on the vertex $v_2$ in Definition \ref{newmov}.
Let
\begin{equation}\label{O}
f_O(\alpha):=\frac{2\alpha}{1-\sqrt{\frac{\alpha}{1-(r+1)\alpha}}}.
\end{equation}
Observe that $f_O(\alpha) \to 0+$ if $\alpha \to 0+$, $f_O(\alpha) \to +\infty$ if $\alpha \to \frac{1}{r+2}-$,
and $f_O(\alpha)$ is strictly increasing in $\alpha \in (0, \frac{1}{r+2})$.
So, there exists a unique $\alpha_1 \in (0, \frac{1}{r+2})$ such that $f_O(\alpha_1)=1$.
Hence $O_{m}$ is consistently $\alpha_1$-normal, and $\rho(O_{m})=\alpha_1^{-\frac{1}{k}}$ by Lemma \ref{comp}.

%We now show $f_P(\alpha_1)>1$, which will imply that $\alpha_0<\alpha_1$.
As $f_O(\alpha_1)=1$, by Eq. (\ref{O}), we have
\begin{equation} \label{ra1}
r \alpha_1=1-\alpha_1-\frac{\alpha_1}{(1-2\alpha_1)^2}.
\end{equation}
Substituting (\ref{ra1}) into Eq. (\ref{P}), we have
\begin{align*}
f_P(\alpha_1)&=\frac{2 \alpha_1}{1-\sqrt{\frac{\alpha_1}{1-\alpha_1}}}+1-\alpha_1-\frac{\alpha_1}{(1-2\alpha_1)^2}\\
&=1+\frac{2 \alpha_1 \left((1-2 \alpha_1)\sqrt{\alpha_1(1-\alpha_1)}-\alpha_1\right)}{(1-2 \alpha_1)^2}
\end{align*}
Let $\psi(\alpha):=(1-2 \alpha)\sqrt{\alpha(1-\alpha)}-\alpha$.
When $\alpha \in (0, \frac{1}{2})$, $\psi(\alpha)>0$ if and only if $$1-\alpha-\frac{\alpha}{(1-2\alpha)^2}>0.$$
By Eq. (\ref{ra1}), surely $\psi(\alpha_1)>0$.
So, $f_P(\alpha_1)>1$, and hence $\alpha_0<\alpha_1$ as $f_O(\alpha)$ and $f_P(\alpha)$ are strictly increasing in $(0,\frac{1}{r+2})$.
By Lemma \ref{comp},
$$ \rho(O_{m})=\alpha_1^{-\frac{1}{k}}<\alpha_0^{-\frac{1}{k}}=\rho(P_{m}).$$
The result follows.\hfill $\square$.

%First consider the root $\alpha_0$ of $f_5(\alpha)-\frac{1}{2}$ when $m=5$.
%Using Mathematica, we have $\alpha_0 \approx 0.205123$.
%
%When $m \ge 6$, $f_m(\frac{1}{2(m-3)}) < \frac{1}{2}$, implying that  $\alpha_0 > \frac{1}{2(m-3)}$.
%
%We now define a  weighted incidence matrix $\bar{B}$ of $T_{m,3}$.
%For each pendent vertex $u$ incident to the unique edge $e$, $\bar{B}(v,e)=1$.
%For each pendent edge $e$ attached at $v \in \{v_2,v_3\}$), define $\bar{B}(v,e)=\alpha_0$.
%Let $\beta_0=\sqrt{\frac{\alpha_0}{1-(m-3)\alpha_0}}$, and $t=1-(m-4)\alpha_0$.
%Define
%$$\bar{B}(v_1,e_1)=\bar{B}(v_1,e_3)=\frac{\alpha_0}{(1-\beta_0)},
%\bar{B}(v_2,e_1)=t(1-\beta_0), \bar{B}(v_2,e_2)=t \beta_0,
%\bar{B}(v_3,e_2)=\beta_0, \bar{B}(v_3,e_3)=1-\beta_0.$$
%Then $\bar{B}$ is consistent.
%To make $T_{m,3}$ is $\alpha_0$-supernormal, we require $t \beta_0^2 \le \alpha$.

\begin{lemma}\label{4cycle}
For $m \ge 4$, $\rho(S^k_{m,4})< \rho(O_{m}).$
\end{lemma}

\emph{Proof.} Label the partial vertices and edges of $S^k_{m,4}$ as in Fig. \ref{8graphs}.
Let $e_1,e_2$ be two non-pendent edges of $S^k_{m,4}$ incident to the vertex of maximum degree,
and let $\{u_1,u_2\} \subseteq e_1$ and  $\{v_1,v_2\} \subseteq e_2$, where $u_2,v_2$ are pendent (cored) and $u_1,v_1$ are non-pendent.
Let $e_3$ be the edge incident to $v_1$ except $e_2$.
Then $e_1,e_2$ satisfy the condition in Definition \ref{newmov}, and by moving $e_3$ from $v_1$ to $u_2$, we get a hypergraph ${S^k_{m,4}}_{e_1,e_2}$ which is isomorphic to $O_{m}$.
By Lemma \ref{Ymov}, we have
$$ \rho(S^k_{m,4})< \rho({S^k_{m,4}}_{e_1,e_2})=\rho(O_{m}).$$
The result follows. \hfill $\square$

%\begin{lemma}\label{Uk}
%For $m \ge 5$, $\rho(U^k_{m,1})< \rho(T^k_{m,1}).$
%\end{lemma}
%
%\emph{Proof.} Label the partial vertices and edges of $U^k_{m,1}$ as in Fig. \ref{8graphs},
%   where $v_1e_1v_2e_2v_3e_3v_1$ is the $3$-cycle, and $w \ne v_2$ and is the vertex to which a pendent edge is attached.
%Let $x$ be a Perron vector of $U^k_{m,1}$.
%If $x_{v_3} \ge x_{w}$, relocating the pendent edge incident to $w$ from $w$ to $v_3$, we will arrive at the hypergraph $T^k_{m,1}$.
%By Lemma \ref{relocate}, $\rho(T^k_{m,1})> \rho(U^k_{m,1})$.
%Otherwise, $x_{v_3} < x_{w}$, moving the edge $e_3$ from $v_3$ to $w$, we also arrive at the hypergraph $T^k_{m,1}$,
%implying that $\rho(T^k_{m,1})> \rho(U^k_{m,1})$ by Lemma \ref{mov}.
%\hfill $\square$

\begin{lemma}\label{Tm2}
For $m \ge 8$, $\rho(T^k_{m,2})< \rho(U^k_{m,1})$.
\end{lemma}

\emph{Proof.} By \cite[Lemma 9]{gsg} or \cite[Theorem 6]{gjm},
if $m \ge 8$, then $\rho(T_{m,2})< \rho(U_{m,1})$.
The result now follows $\rho(T^k_{m,2})=\rho(T_{m,2})^{\frac{2}{k}}<\rho(U_{m,1})^{\frac{2}{k}}= \rho(U^k_{m,1})$ by Lemma \ref{pow}.
\hfill $\square$

\begin{lemma}\label{UQ}
For $m \ge 5$, $\rho(U^k_{m,1})< \rho(Q_{m})$.
\end{lemma}

\emph{Proof.} 
Label the partial vertices and edges of $U^k_{m,1}$ as in Fig. \ref{8graphs}.
Let $e_4$ be the pendent edge incident to $w$, and let $e_5$ be the non-pendent edge incident to $w$.
Now $e_2,e_5$ satisfy the condition in Definition \ref{newmov}.
Let $u$ be a pendent (cored) vertex of $e_2$.
Moving $e_4$ from $w$ to $u$, we get a hypergraph ${U^k_{m,1}}_{e_2,e_5}$ isomorphic to $Q_m$.
By Lemma \ref{Ymov}, we have
$$ \rho(U^k_{m,1})< \rho({U^k_{m,1}}_{e_2,e_5})=\rho(Q_{m}).$$
The result follows. \hfill $\square$

%   where $v_1e_1v_2e_2v_3e_3v_1$ is the $3$-cycle, and $w$ is the vertex of the edge $e_3$ to which a hyperstar is attached.
%We define a weighted incidence matrix $B$ of $U^k_{m,1}$ as follows.
%For each pendent vertex $v$ incident to the unique edge $e$, $B(v,e)=1$.
%%Let $w$ be the vertex to which all pendent edges are attached.
%For each pendent edge $e$ attached at $v \in \{v_2,w\}$, define $B(v,e)=\alpha$.
%Let $r:=m-4 \ge 1$, and let $e_4$ be the non-pendent edge incident to $w$.
%Define
%\begin{align*}
%B(w,e_4)&=1-\alpha, B(v_2,e_4)=\frac{\alpha}{1-\alpha}, B(v_2,e_1)=B(v_2,e_2)=\beta, \\
%B(v_1,e_1)&=B(v_3,e_2)=\frac{\alpha}{\beta}, B(v_1,e_3)=B(V_3,e_3)=1-\frac{\alpha}{\beta}.\\
%\end{align*}
%It is easily seen $B$ is consistent.
%To make $O_{m}$ be $\alpha$-normal, we require

We now determine the linear unicyclic hypergraph with the second  or third largest spectral radius among all linear unicyclic hypergraph with $m$ edges.
%When $m=3$, there is only one linear unicyclic hypergraph, i.e. a linear $3$-cycle.
We need the following result.

\begin{lemma}{\em \cite{FanTPL}}\label{1st}
$(1)$ Among all unicyclic linear $k$-uniform hypergraphs with
$m \ge 4$ edges and girth $g$, the power hypergraph $S^k_{m,g}$ is the unique maximizing
hypergraph.

$(2)$ For $g \ge4$, $\rho(S^k_{m,g})<\rho(S^k_{m,g-1})$.

$(3)$ Among all unicyclic linear $k$-uniform hypergraphs with $m \ge 4$ edges, $S^k_{m,3}$ is the unique maximizing hypergraph.
\end{lemma}

%Let $G$ be a linear unicyclic $k$-uniform hypergraphs with $m > 3$ edges.
%Then $G$ is obtained from a linear cycle $C$ by attached some hypertrees at its vertices.
%If a hypertree $T$ is attached at a cored vertex of $C$, then $T$ is called a \emph{cored branch};
%otherwise, $T$ is called a \emph{non-cored branch}.

\begin{theorem}\label{2nd}
Among all linear unicyclic $k$-uniform hypergraphs with $m \ge 5$ edges, $T^k_{m,1}$ is the unique hypergraph with the second largest spectral radius.
\end{theorem}

\emph{Proof.}
Let $G$ be a hypergraph with the second largest spectral radius among all linear unicyclic $k$-uniform hypergraphs with $m > 3$ edges.
Surely $G \ne S^k_{m,3}$ by Lemma \ref{1st}.
In the following we call a hypergraph \emph{proper} if it is not equal to $S^k_{m,3}$.
Suppose $G$ has girth $g$.
We assert that $g=3$; otherwise, by Lemma \ref{1st} and Lemma \ref{4cycle},
\begin{equation}\label{g4}
\rho(G) \le \rho(S^k_{m,g}) \le \rho(S^k_{m,4})<\rho(O_{m}).
\end{equation}
So $G$ is obtained from a cycle $C$ of length $3$ by attached some hypertrees at its vertices.

Let $x$ be a Perron vector of $G$. The result will follows by the following cases.

{\bf Case 1.} Exactly one hypertree is attached to some vertex of $C$.
Let $T_u$ be such hypertree attached at $u$ of $C$.
Write $G=C(u) \ast T_u(u)$.

{\bf Case 1.1.} $u$ is a cored vertex of $C$.
Then $u$ is the unique vertex of $T_u$ such that $x_u=\max\{x_v: v \in V(T_u)\}$.
Otherwise, let $v \in V(T_u) \backslash \{u\}$ such that $x_v \ge x_u$.
Relocating $C$ from $u$ to $v$, we will get a proper hypergraph $G'$ which holds $\rho(G') > \rho(G)$ by Lemma \ref{relocate}, a contradiction.

We assert $T_u$ is a hyperstar with center $u$.
Otherwise there exists a pendent edge $e$ of $T_u$ incident to a non-cored vertex $w \ne u$.
Relocating the edge $e$ from $w$ to $u$, we also get a proper hypergraph but with a larger spectral radius, a contradiction.
So $G=O_{m}$ in this case.
However, by Lemma \ref{OP}, $\rho(O_{m})< \rho(P_{m})$.
So this case cannot happen.

{\bf Case 1.2.} $u$ is a vertex of $C$ of degree two.
Then $u$ is the unique vertex of $T_u$ such that $x_u=\max\{x_v: v \in V(T_u)\}$.
Otherwise, let $v \in V(T_u) \backslash \{u\}$ such that $x_v \ge x_u$, and let $e$ be an edge of $C$ incident to $u$.
Moving $e$ from $u$ to $v$, we also get a proper hypergraph $G'$ with girth at least $4$ and a larger spectral radius, a contradiction.
We assert $G=U^k_{m,1}$.
Otherwise, as $G \notin \{S^k_{m,3},U^k_{m,1}\}$, there exists a pendent edge $e$ of $T_u$ incident to a non-cored vertex $w \ne u$.
Relocating the edge $e$ from $w$ to $u$, we also get a proper hypergraph but with a larger spectral radius, a contradiction.
However, by Lemma \ref{UQ}, $\rho(U^k_{m,1}) < \rho(Q_m)$.
So this case cannot happen.

{\bf Case 2.} At least two hypertrees are attached at different vertices of $C$.
Suppose that there are $s$ vertices, say $v_1,\ldots,v_s$ of $C$, are attached $s$ hypertrees, where $s \ge 3$.
Without loss of generality, assume that $x_{v_1}=\max\{x_{v_i}: i=1,\ldots,s\}$.
Relocating the hypertree attached at $v_s$ from $v_s$ to $v_1$, we will get a proper hypergraph with larger spectral radius, a contradiction.
So, there are exactly two hypertrees, say $T_u$ and $T_v$, attached at $u$ and $v$ of $C$ respectively.
By a similar discussion as in Case 1.1, $T_u$ and $T_v$ are hyperstars with center $u$ and $v$ respectively.

{\bf Case 2.1.} Both $u$ and $v$ are cored vertices of $C$.
Without loss of generality, $x_u \ge x_v$.
Relocating $T_v$ from $v$ to $u$, we will get the hypergraph $O_{m}$ with larger spectral radius, implying this case also cannot happen.

{\bf Case 2.2.} One of $u,v$ is a cored vertex and the other is non-cored vertex of $C$.
Without loss of generality, $u$ is cored and $v$ is non-cored.
First assume that $u, v$ are adjacent. If $x_u \ge x_v$, relocating $T_v$ from $v$ to $u$, we will arrive the hypergraph $O_{m}$, a contradiction.
%However, $\rho(O_{m,3})< \rho(P_{m,3})$, implying this case cannot happen.
Otherwise, $x_u < x_v$, moving all pendent edges except one arbitrarily specified edge incident with $u$ from $u$ to $v$
if there exists more than one pendent edges incident to $u$,
we will get the hypergraph $Q_{m}$ with larger spectral radius. So $G=Q_m$ in this case.
However, by Lemma \ref{QT}, $\rho(Q_{m})< \rho(T^k_{m,1})$, implying this case also cannot happen.

Secondly assume that $u,v$ are not adjacent.
Similarly, if $x_u \ge x_v$, relocating $T_v$ from $v$ to $u$, we will arrive the hypergraph $O_{m}$, a contradiction.
%implying this case cannot happen as $\rho(O_{m,3})< \rho(P_{m,3})$.
If $x_u < x_v$, moving all pendent edges except one arbitrarily specified edge incident with $u$ from $u$ to $v$
if there exists more than one pendent edges incident to $u$,
we will get the hypergraph $P_{m}$ with larger spectral radius. So $G=P_m$ in this case.
However, by Lemma \ref{PQ}, $\rho(P_{m})< \rho(Q_{m})$, implying this case also cannot happen.

{\bf Case 2.3.} Both $u,v$ are non-cored vertices of $C$.
Without loss of generality, assume that $x_u \ge x_v$.
Moving all pendent edges incident with $v$ except one arbitrarily specified edge from $v$ to $u$ (if there exists more than one pendent edges incident to $v$),
we will get the hypergraph $T^k_{m,1}$ with larger spectral radius. So $G=T^k_{m,1}$ in this case.
\hfill $\square$

\begin{theorem}\label{3rd}
Among all linear unicyclic $k$-uniform hypergraphs with $m \ge 8$ edges, $Q_m$ is the unique hypergraph with the third largest spectral radius.
\end{theorem}

\emph{Proof.}
Let $G$ be a hypergraph with the third largest spectral radius among all linear unicyclic $k$-uniform hypergraphs with $m \ge 8$ edges.
Surely $G \notin \{S^k_{m,3}, T^k_{m,1}\}$ by Lemma \ref{1st} and Lemma \ref{2nd}.
In the following we call a hypergraph \emph{proper} if it is not equal to $S^k_{m,3}$ or $T^k_{m,1}$.
By (\ref{g4}), we know $G$ has girth $3$, and $G$ is obtained from a cycle $C$ of length $3$ by attached some hypertrees at its vertices.
Let $x$ be a Perron vector of $G$.
We now follow the routine of the proof of Theorem \ref{2nd}.

{\bf Case 1.} Exactly one hypertree is attached to some vertex, say $u$, of $C$.
We assert $u$ is a non-cored vertex, and hence $G=U^k_{m,1}$ by Case 1.2 in the proof of Theorem \ref{2nd}.
However, by Lemma \ref{UQ}, $\rho(U^k_{m,1}) < \rho(Q_m)$, a contradiction.
Otherwise, if $u$ is a cored vertex, then from Case 1.1 of Theorem \ref{2nd}, $G=O_{m}$ and $\rho(O_{m})< \rho(P_{m})$ by Lemma \ref{OP}, a contradiction.

{\bf Case 2.} At least two hypertrees are attached at different vertices of $C$.
By a similar discussion as in Case 2 of Theorem \ref{2nd}, $G$ is obtained from $C$ by attaching two hyperstars $T_u, T_v$ at $u,v$ of $C$,
where $u,v$ are the centers of $T_u, T_v$ respectively.

{\bf Case 2.1. } If $u,v$ are both cored vertices of $C$, then $\rho(G) < \rho(O_{m})$, a contradiction.

{\bf Case 2.2.} Suppose that $u$ is cored and $v$ is non-cored.
If $u,v$ are adjacent, then $\rho(G) < \rho(O_{m})$ or $G=Q_{m}$ depending on whether $x_u \ge x_v$  or not.
So $G=Q_{m}$ in this case.
If $u,v$ are not adjacent, then $\rho(G) < \rho(O_{m})$ or $G=P_{m}$ depending on whether $x_u \ge x_v$  or not.
However, $\rho(P_{m})< \rho(Q_{m})$ by Lemma \ref{PQ}, a contradiction.

{\bf Case 2.3.} Finally suppose that $u,v$ are non-cored vertices of $C$.
As $G \notin \{ S^k_{m,3}, T^k_{m,1}\}$, both $T_u, T_v$ have at least $2$ pendent edges.
Without loss of generality, assume that $x_u \ge x_v$.
Moving all pendent edges incident with $v$ except two arbitrarily specified edge from $v$ to $u$ if there exist more than two pendent edges incident to $v$,
we will get the hypergraph $T^k_{m,2}$ with a larger spectral radius.
So $G=T^k_{m,2}$.
So $\rho(T^k_{m,2}) < \rho(U^k_{m,1})$ by Lemma \ref{Tm2}, a contradiction.

The result now follows by the above discussion.
\hfill $\square$

\end{document}